\documentclass{article}
\usepackage{amsmath,amsthm,amssymb}
\usepackage{epsfig}

\begin{document}
\newtheorem{thm}{Theorem}
\newtheorem{lem}{Lemma}
\newtheorem{cor}{Corollary}
\theoremstyle{definition}
\newtheorem{nota}{Notation}
\newtheorem{rem}{Remark}

\title{A Sharp Condition for the Loewner Equation to Generate Slits}

\author{Joan R. Lind}

\date{}

\maketitle

\begin{abstract}

D. Marshall and S. Rohde have recently shown that there exists $C_0 >0$ so 
that the Loewner equation generates slits whenever the driving term 
is H\"older continuous with exponent $\frac{1}{2}$ and norm less than 
$C_0$ \cite{RM}.  
In this paper, we show that the maximal value for $C_0$ is 4.

\end{abstract}

\subsection*{Introduction}

\newcommand{\norm}[1]{\lVert#1\rVert}
\newcommand{\abs}[1]{\lvert#1\rvert}

When Loewner introduced his namesake differential equation in 1923, it
greatly impacted the theory of univalent functions.  A univalent function $f$ 
is a
conformal map of the unit disc, normalized by $f(0)=0$ and $f'(0)=1$.  In
other words, it has the following power series representation in the unit 
disc:

$$f(z)=z+a_2z^2+a_3z^3+ \cdots.$$ 
In 1916 Bieberbach \cite{BB} had shown that $\abs{a_2} \leq 2$ and had
conjectured that $\abs{a_n} \leq n$ for all $n$.  It was Loewner's
differential equation that led to a proof of the case $n=3$ in 1923.  
See \cite{A1} or \cite{D} for a proof of this and for more classical 
applications of the Loewner equation.  When the
Bieberbach conjecture finally was proved for general $n$ in 1985 by de 
Branges \cite{dB}, the Loewner equation again played a key role.

In addition to its importance in the theory of univalent functions, the 
Loewner 
differential equation has gained recent prominence with the 
introduction of a stochastic process called ``Stochastic Loewner 
Evolution", or SLE, by O. Schramm \cite{SLEintro}.   Many results in this 
fast-growing 
field can be found in the recent work of mathematicians such as
Lawler, Rohde, Schramm, Smirnov, and Werner.  See \cite{Kad1} for a survey 
paper with an extensive bibliography.

In the next two sections, we will introduce two formulations of the
deterministic Loewner differential equation, the halfplane version and
the disc version.  This is followed by a discussion of some problems 
associated with the
geometry of the solutions to the Loewner equation.
The rest of the paper is concerned with proving Theorem \ref{thm:1} 
below, which builds upon D. Marshall and S. Rohde's recent work \cite{RM} 
concerning when the Loewner equation can generate slits.  
The fifth section contains examples and lemmas
related to a natural obstacle to generating slits, the sixth 
section includes lemmas
about conformal welding and the Loewner equation, and the 
final section is the proof of Theorem \ref{thm:2}, which 
is equivalent to Theorem \ref{thm:1}.

\subsection*{The Loewner equation in the halfplane}

Let $\gamma(t)$ be a simple continuous curve in 
$\mathbb{H} \cup \{0\}$ 
with $\gamma(0) =0$ and $ t \in [0,T]$.  
Then there is a unique conformal map $g_t: \mathbb{H} \setminus 
\gamma[0,t] 
\rightarrow \mathbb{H}$ with the following normalization, called the 
hydrodynamic normalization, near infinity:

$$g_t(z)=z+\frac{c(t)}{z} + O\left(\frac{1}{z^2}\right).$$
It is an easy exercise to check that $c(t)$ is continuously increasing in 
$t$ and that 
$c(0) = 0$.  Therefore $\gamma$ can be reparametrized so that $c(t) = 2t$.  
Assuming this normalization, one can show that $g_t$ satisfies the 
following form of 
Loewner's differential equation:
for all $ t \in
[0,T]$  and all $z \in \mathbb{H} \setminus \gamma[0,t]$, 

$$ \frac{\partial}{\partial t} g_t(z) = \frac{2}{g_t(z)- \lambda(t)},$$

$$g_0(z) = z,$$
where $\lambda$ is a continuous, real-valued function.  Further, it can be 
shown that $g_t$ extends continuously to $\gamma(t)$  and 
$g_t(\gamma(t))$ equals $\lambda(t)$.

On the other hand, if we start with a continuous $\lambda :[0,T] 
\rightarrow \mathbb{R}$, we can consider the following initial value 
problem for each $z \in \mathbb{H}$:

\begin{gather}\label{bw}
\frac{\partial}{\partial t} g(t,z) = \frac{2}{g(t,z)- \lambda(t)},\\
g(0,z) = z. \notag
\end{gather}
For each $z \in \mathbb{H}$ there is some time 
interval $[0, s)$ for which a solution $g(t, z)$ exists.  Let $T_z= 
\sup\{s \in [0,T]: g(t,z)$ exists on $[0,s)\}$.  Set $G_t = 
\{z \in \mathbb{H} : T_z > t\}$ and $g_t (z) = g(t,z)$.  Then one can 
prove that the set $G_t$ is a simply connected subdomain of 
$\mathbb{H}$ and $g_t$ is the unique conformal map from $G_t$ onto 
$\mathbb{H}$ with the following normalization near infinity:

$$g_t(z)=z+\frac{2t}{z} + O\left(\frac{1}{z^2}\right).$$
The function $\lambda(t)$ is called the driving term, and the domains 
$G_t$ as well as the functions $g_t$ are said to be generated by 
$\lambda$.  

The domains $G_t$ generated by a continuous driving term $\lambda$ are not
necessarily slit-halfplanes, i.e. domains of the form $\mathbb{H}
\setminus \gamma[0,t]$, for some simple continuous curve $\gamma$ in
$\mathbb{H} \cup \{\gamma(0)\}$ with $\gamma(0) \in \mathbb{R}$.  We 
will give an 
example later in 
the paper where
a non-slit-halfplane is generated by a driving term which is not only
continuous but also is in Lip$(\frac{1}{2})$.  Recall that
Lip$(\frac{1}{2})$ is the space of H\"older continuous functions with
exponent $\frac{1}{2}$, that is the space of functions $\lambda(t)$
satisfying $\abs{\lambda(s)-\lambda(t)} \leq c \abs{s-t}^{1/2}$,
with $ \norm{\lambda}_{\frac{1}{2}}$ denoting the smallest such $c$. 
The necessary and sufficient
condition for a decreasing family of domains $\{G_t\}$ to be generated by
a continuous driving term can be found in Section 2.3 of \cite{LSW11}.

\subsection*{The Loewner equation in the disc}

The setup for the disc version of the Loewner equation is similar to that 
of the halfplane version, but the normalization will be at an 
interior point rather than at a boundary point.
For the unit disc $\mathbb{D}$ slit by a simple curve $\gamma(t)$ in 
$\mathbb{D} \cup \{1 \}$ with
$\gamma(0)=1$ and $\gamma(t) \neq 0$ for any $t$,    there 
is a unique 
family 
of conformal 
maps 
$\{g_t\}$ so that $g_t: \mathbb{D} \setminus \gamma[0,t]
\rightarrow \mathbb{D}$ with the normalizations $g_t(0)=0$ and 
$g_t'(0)>0$.  Further, by reparametrizing $\gamma$ if necessary, we can 
assume that $g_t'(0)=e^t$.  If we again set $\lambda(t)=g_t(\gamma(t))$, 
then 

\begin{gather}\label{deq}
\frac{\partial}{\partial t} g_t(z) = g_t(z) 
\frac{\lambda(t)+g_t(z)}{\lambda(t)-g_t(z)},\\
g_0(z)=z. \notag
\end{gather}

Given any continuous function 
$\lambda: [0, T] \rightarrow \partial \mathbb{D}$, we can solve the 
initial 
value problem \eqref{deq} for $z \in \mathbb{D}$.   As in the halfplane 
version, this will 
generate a family of 
conformal maps $\{g_t\}$ which map from a  simply 
connected 
subdomain of the unit disc onto the unit disc and which are normalized by 
$g_t(0)=0$ and $g_t'(0)=e^t$.

\subsection*{Some results}

We return to the halfplane version of the Loewner equation, which will be 
the setting for the rest of this paper.
For $\kappa
\geq 0$, set $\lambda(t) = \sqrt\kappa B_t$, where $B_t$ is standard
Brownian motion.  Then chordal SLE$_{\kappa}$ is the random family of
conformal maps generated by $\lambda$, that is, the family
of maps solving the
following stochastic differential equation:
 
\begin{gather*}
\frac{\partial}{\partial t} g_t(z) = \frac{2}{g_t(z)- \sqrt\kappa B_t},\\
g_0(z) = z.
\end{gather*}

For SLE, it is possible to define an almost surely continuous path
$\gamma: [0,\infty) \rightarrow \overline{\mathbb{H}}$ such that the
domains $G_t$ generated by
$\lambda(t)
= \sqrt\kappa B_t$ are the unbounded components of $\mathbb{H} \setminus
\gamma[0,t]$ for every $t \geq 0$.  See \cite{RS} and, for the case
$\kappa=8$, \cite{LSW14}.  Further, S. Rohde and O. Schramm \cite{RS} have
shown the following classification:

\begin{enumerate}
\item
For $\kappa \in [0,4]$, $\gamma (t)$ is almost surely a simple
path contained in $\mathbb{H} \cup \{0\}$.
\item
For $\kappa \in (4,8)$, $\gamma (t)$ is almost surely a
non-simple path.
\item
For $\kappa \in [8,\infty)$, $\gamma(t)$ is almost surely a
space-filling curve.
\end{enumerate}

This result motivates a question in the deterministic setting.
Can we classify 
the kinds of domains generated by a driving term $\lambda$ 
in terms of 
some characteristic of $\lambda$?  There is 
only a partial understanding 
of this question.  In the case of a domain slit by an analytic slit, the 
driving term is real analytic, and if the slit is 
$C^n$ then the driving term is at least $C^{n-1}$.  See \cite{EE} and 
\cite{BLW}.

D. Marshall and S. Rohde address the question of when the 
generated 
domains $G_t$ are quasislit-halfplanes in \cite{RM}, where a 
quasislit-halfplane is the image of $ \mathbb{H} \setminus [0,i]$
under a quasiconformal mapping fixing $\mathbb{H}$ and $\infty$.
They prove the following:

\begin{thm}\label{rmthm}
If $G_t$ is a quasislit-halfplane for all $t$, then $\lambda  \in$
Lip$(\frac{1}{2})$.  Conversely,
there exists $C_0$ such 
that if the driving term 
$\lambda  \in$ Lip$(\frac{1}{2})$ with $ \norm{\lambda}_{\frac{1}{2}} < 
C_0$, then $G_t$ is a quasislit-halfplane for all $t$.
\end{thm}
 Although they work 
with the  
technically more challenging disc version of the Loewner equation, their 
techniques carry over to prove the result in the halfplane version as 
well.  In the remainder of this paper, working with the halfplane version 
of Loewner's 
equation, we will show that the maximal value for $C_0$ is 4.  

\begin{thm}\label{thm:1}
If $\lambda \in$ Lip$(\frac{1}{2})$ with $ \norm{\lambda}_{\frac{1}{2}} < 
4$, 
then the domains
$G_t$ generated by $\lambda$ are quasislit-halfplanes.  
\end{thm}

Further, for each $c \geq 4$,
there exists a driving term
$\lambda \in$ 
Lip$(\frac{1}{2})$ with $ \norm{\lambda}_{\frac{1}{2}} =
c$ so that $\lambda$ does not generate slit-halfplanes.  We will see 
examples of this in the next section.  Similar examples were discovered 
independently by L. Kadanoff, W. Kager, and B. Nienhuis \cite{Kad2}.  
Their work also includes descriptions and pictures of the generated 
domains.

There is another version of the Loewner equation in the halfplane.  Let 
$\xi : [0,T] \rightarrow \mathbb{R}$ be continuous and consider 
the following initial value problem, in which a negative 
sign has been introduced on the righthand side of \eqref{bw}:

\begin{gather}\label{fw}
\frac{\partial}{\partial t} f(t,z) = \frac{-2}{f(t,z)-
\xi(t)},\\
f(0,z) = z \notag
\end{gather}
for $z \in \mathbb{H}$.  In this case, for each $z \in 
\mathbb{H}$, the solution $f(t,z)$ exists for all $t \in [0,T]$.  
Setting $f_t(z) = f(t,z)$, we have that $f_t$ is defined on all of 
$\mathbb{H}$.  As in the previous case, it can be shown that $f_t$ is a 
conformal map from $\mathbb{H}$ into $\mathbb{H}$, and near infinity it 
has the form

$$f_t(z)=z+\frac{-2t}{z} + O(\frac{1}{z^2}).$$
We think of the funtions $f_t$ as
being generated by ``running
time backwards."

These two forms of Loewner's differential equation are related.  Given  a
continuous function $\lambda$ on $[0,T]$, set $\xi(t)= \lambda(T-t)$.  Let 
$g_t$ be the functions generated by $\lambda$ from \eqref{bw}, and let 
$f_t$ be 
the functions generated by $\xi$ from \eqref{fw}.  It is not true that 
$f_t(z)=g_t^{-1}(z)$ for all $t \in [0,T]$,  but it is true that 
$f_T(z)=g_T^{-1}(z)$.  
Therefore Theorem \ref{thm:1} is 
equivalent to the following:

\begin{thm}\label{thm:2}
If $\xi \in$ Lip$(\frac{1}{2})$ with $ \norm{\xi}_{\frac{1}{2}} <
4$,
then
$f_t(\mathbb{H})$ is a quasislit-halfplane for all t, where $f_t$ are 
the maps generated by $\xi$.
\end{thm}

\subsection*{When the singularity catches solutions}

Let  $\lambda  \in$ Lip$(\frac{1}{2})$ and suppose that 
the domains $G_t$ generated by $\lambda$ are slit-halfplanes.  Then the 
maps $g_t$ extend continuously to $\mathbb{R} \setminus \{\lambda(0)\}$.  
Thus for each $x_0 \in \mathbb{R} \setminus \{\lambda(0)\}$, $x(t) := 
g_t(x_0)$ 
is a solution to the following real-valued initial value problem:

\begin{gather}\label{bwr}
\frac{\partial}{\partial t} x(t) = \frac{2}{x(t)- \lambda(t)},\\
x(0) = x_0. \notag
\end{gather}
Further, if $\lambda$ is defined on $[0,T]$, then $x(t) 
\neq \lambda(t)$ 
for any $t \in [0,T]$, since otherwise \eqref{bwr} would fail to have a 
solution 
for all $t \in [0,T]$.  

Note that if $x_0 > \lambda(0)$, then $ 
\frac{\partial}{\partial t} x(t) >0$ as long as $x(t) \neq \lambda(t)$.  
So two things can happen: either $x(t)$ continues to move to the right, 
staying strictly larger than the driving term, or the driving term 
moves fast enough to ``catch" $x(t)$ and there is some time $t_0$ 
where $x(t_0)=\lambda(t_0)$.  The case when $x_0 < \lambda(0)$ is similar 
but with 
$x(t)$ moving to the left.  Thus, when the domains generated are 
slit-halfplanes, we see that $\lambda(t)$ cannot ``catch"  
any solution $x(t)$ to \eqref{bwr}.  

To build our intuition, let us briefly consider a particular example.  Let 
$G_t= \mathbb{H} \setminus \gamma[0,t]$, where $\gamma$ parametrizes  
the upper half-circle of radius 
$\frac{1}{2}$ centered at $\frac{1}{2}$, as pictured in Figure 1.  
In this 
case it is possible, although unpleasant, to compute the maps $g_t$ and to 
ascertain that the driving term generating this scenario is the function 
$\lambda(t)=\frac{3}{2}-\frac{3}{2}\sqrt{1-8t}$, for $t \in [0, 
\frac{1}{8}]$.  The time $t=\frac{1}{8}$ corresponds to the moment that 
the 
circular arc touches back on the real line, and $G_{\frac{1}{8}} = 
\mathbb{H} \setminus D(\frac{1}{2}, \frac{1}{2})$.

\begin{figure}
\centering
\epsfig{file=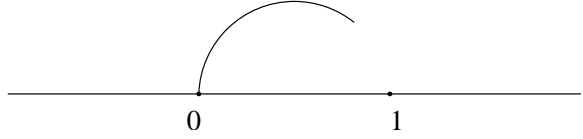}
\caption{One of the domains generated by
$\lambda(t)=\frac{3}{2}-\frac{3}{2}\sqrt{1-8t}$  }
\end{figure}

For $t \in [0, \frac{1}{8}-\epsilon]$, the domains $G_t$ are 
slit-halfplanes, and therefore for any $x_0 \neq 0$, the solutions $x(t)$ 
to \eqref{bwr} exist on this time interval.  What happens to these 
solutions when $t=\frac{1}{8}$?  Clearly, 
$g_{\frac{1}{8}}$ extends only to $\mathbb{R} \setminus [0,1]$.  That is, 
on $[0,\frac{1}{8}]$, solutions to \eqref{bwr} exist only for $x_0 >1$ or 
$x_0 <0$.  So if $x_0 \in (0,1]$, the function $x(t)$ resulting from 
\eqref{bwr} must be ``caught" by $\lambda$ at time $t=\frac{1}{8}$.  For 
example,  it is easy to check that the solution to \eqref{bwr} when 
$x_0=1$ is $x(t)=\frac{3}{2}-\frac{1}{2}\sqrt{1-8t}$.  Here we see that 
$x(\frac{1}{8})=\frac{3}{2}=\lambda(\frac{1}{8})$.

To determine an upper bound on the constant $C_0$ in Theorem \ref{rmthm},
we can analyze the situations where this ``catching" could occur, since
this implies that the family of domains $G_t$ is not a family of
slit-halfplanes.  In the example above,
$\norm{\lambda}_{\frac{1}{2}}=3\sqrt2$, which indicates that $C_0 \leq
3\sqrt2$.  Moreover, for any $c \geq 4$, it is easy to give an example of
a driving term $\lambda$ with $\norm{\lambda}_{\frac{1}{2}} = c$ so that
$\lambda$ can "catch" a function $x(t)$ generated by \eqref{bwr}
for some $x_0$.  Let $\lambda(t)= c-c\sqrt{1-t}$ and $x(t)=c-a\sqrt{1-t}$
where $a=\frac{1}{2} (c+\sqrt{c^2-16})$.  In particular, when $c=4$, then
$\lambda(t)=4-4\sqrt{1-t}$ and $x(t)=4-2\sqrt{1-t}$.  One can check that
$x(t)$ is a solution to \eqref{bwr} with $x_0=c-a>0$.  However
$x(1)=c=\lambda(1)$.  Therefore, since $\lambda(t)$ has ``caught" $x(t)$,
$\lambda$ cannot generate slit-halfplanes.  This implies that the constant
$C_0$ in Theorem \ref{rmthm} cannot be greater than 4.

In contrast to the examples above, the following lemma shows that if 
$\lambda$ can "catch" some $x(t)$, then $
\norm{\lambda}_{\frac{1}{2}}
\geq 4$.  To make things slightly simplier, we take advantage of the fact 
that the halfplane version of the Loewner equation 
satisfies a useful scaling property: If $\lambda (t)$ and $x(t)$ 
satisfy equation \eqref{bwr}, then $\hat{\lambda}(t) := \frac{1}{r} 
\lambda(r^2t)$ and $\hat{x}(t) := \frac{1}{r} x(r^2t)$ also satisfy 
equation \eqref{bwr}.  Verifying this is an easy exercise.  
Using this 
scaling property, we can assume that if a "catching" occurs, 
then it happens at time 1.  More precisely, if
$x(t_0) = \lambda(t_0)$ and $x(t) \neq \lambda (t)$ for $t < 
t_0$, then without loss of generality $t_0=1$.  Also, nothing is lost by 
assuming  that $\lambda(0)=0$ and $x_0>0$.

\begin{lem}\label{catchpoints}
Let  $\lambda  \in$ Lip$(\frac{1}{2})$ with $\lambda(0)=0$ and let $x_0
>0$.  Suppose that $x(t)$ is a solution to
\eqref{bwr} and 
that
$x(1)=\lambda(1)$.  Then  $ \norm{\lambda}_{\frac{1}{2}} \geq
4$.

\end{lem}

\begin{proof}

Let c=$ \norm{\lambda}_{\frac{1}{2}}$.  
From \eqref{bwr}, we have that $x(t)$ is increasing in $t$.  
So then since $\lambda \in $ Lip$(\frac{1}{2})$,
$$x(t)-\lambda(t) \leq x(1) - \lambda(1) + c \sqrt{1-t} \leq 
c\sqrt{1-t}.$$
From \eqref{bwr} we have
$$ \dot{x}(t) \geq \frac{2}{c\sqrt{1-t}}.$$
Integrating gives that
$$x(1)-x(t) \geq \frac{4}{c} \sqrt{1-t}.$$
Letting $t=0$ and using that $x(1)-x_0<c$, we see that $c - \frac{4}{c} > 
0$ and so $c>2$.  But we 
also have a better estimate for $x(t)$:
$$x(t) \leq x(1) - \frac{4}{c} \sqrt{1-t}.$$
Now using this estimate, we can repeat the above argument.  So
$$x(t)-\lambda(t) \leq (c-\frac{4}{c})\sqrt{1-t},$$
which leads to a new estimate for $\dot{x}(t)$.  Then by integration,
$$x(1)-x(t) \geq \frac{4}{c-\frac{4}{c}} \sqrt{1-t}.$$
This implies that $c-\frac{4}{c-\frac{4}{c}} >0$ and so $c>2\sqrt2$.
Again we also get an improved estimate for $x(t)$:
$$x(t) \leq x(1) -\frac{4}{c-\frac{4}{c}} \sqrt{1-t}.$$

Repeating this procedure $n$ times gives that $h_n(c) > 0$ where $h_n$ is 
recursively defined as follows:
$$h_1(x) = x-\frac{4}{x},$$
$$h_n(x) = x-\frac{4}{h_{n-1}(x)}.$$
Note that $h_1(x)$ is an increasing function from $(0, \infty)$ onto 
$\mathbb{R}$.  It is easy to show inductively that we
can define an increasing sequence $\{x_n\}$ so that $h_n(x_n)=0$, 
and $h_{n+1}(x)$ is an increasing function from $(x_n, \infty)$ onto
$\mathbb{R}$.  Note that we have shown that $x_1=2$ and 
$x_2=2\sqrt 2$.
Since $h_n(c) > 0$ for all $n$, $c>x_n$ for all $n$.  It simply remains to 
show that $x_n \nearrow 4$.

An easy 
inductive argument gives that $h_n(4) \geq 2$ for all $n$.  If $4 \in 
(x_{k-1}, x_k]$ 
for some $k$, then $h_k(4) \leq 0$.  Therefore, the 
increasing sequence 
$\{x_n\}$ is bounded above by 4, and hence there exists some $a \leq 4$ 
such 
that $x_n \nearrow a$.  Now, $h_n(a) > h_n(x_n)=0$ for all $n$.
If $h_k(a) \leq 1$ for some $k$, then $h_{k+1}(a)=a-\frac{4}{h_n(a)} \leq 
0$.  So we must have $h_n(a)>1$ for all $n$.  Since $h_n(a)$ is decreasing 
in $n$ and 
bounded 
below by 1, $h_n(a) \searrow L$ for some $L \geq 1$.  So then,
$$L=\lim_{n \rightarrow \infty} h_n(a) = \lim_{n \rightarrow \infty} 
a-\frac{4}{h_{n-1}(a)}=a-\frac{4}{L}.$$
Solving the above for $L$ gives that
$$L=\frac{a \pm \sqrt{a^2-16}}{2}.$$
Since we know the real-valued limit $L$ exists, we must have $a\geq 4$.  
Hence, $a=4$, completing the proof.

\end{proof}

Note that in the proof above, we have also shown the following: if 
$h_n(c)>0$ for all $n$, then $c 
\geq 4$.  This follows since $h_n(c)>0$ for all $n$ implies that $c>x_n$ 
for all $n$ and since $x_n 
\nearrow 4$ .
We mention this here, since we will use this fact in the proof of the 
next lemma.

Although Lemma \ref{catchpoints} certainly suggests that the maximal 
value for 
$C_0$ is 4, it is not a proof of Theorem \ref{thm:1}.  In theory, there 
may be 
more 
obstacles to generating quasislit-halfplanes than that of the driving term 
catching up to some solution to \eqref{bwr}.   However, we will see that 
this is 
basically the only obstacle.  Refining the above 
argument gives Lemma \ref{catchpoints2}, which combined with techniques in 
\cite{RM} will 
lead to the proof of Theorem \ref{thm:1}.  The idea of Lemma \ref{catchpoints2}
 is that if 
$\lambda$ can get close to catching some $x(t)$, then  $ 
\norm{\lambda}_{\frac{1}{2}}$ must be close to being greater than or 
equal to 4.

\begin{lem}\label{catchpoints2}
Let  $\lambda  \in$ Lip$(\frac{1}{2})$ with $\lambda(0)=0$ and 
$\norm{\lambda}_{\frac{1}{2}} < 4$.  Then there 
exists $\epsilon = \epsilon (\norm{\lambda}_{\frac{1}{2}})> 0$ so that 
$x(1)-\lambda(1) > \epsilon$, where $x(t)$ is the solution to 
\eqref{bwr} with $x_0>0$.

\end{lem}

\begin{proof}

Suppose $x(t)$ is a solution to \eqref{bwr} for some $x_0>0$ so that  
$x(1)-\lambda(1) \leq \epsilon$.  We will show that there exists some 
$\epsilon >0$ so that this leads to a contradiction.  Again, let $c=
\norm{\lambda}_{\frac{1}{2}}$.
As in the previous proof, define $h_n$ recursively by
$$h_1(c) = c- \frac{4}{c},$$
$$h_n(c) = c-\frac{4}{h_{n-1}(c)}.$$
Since $c<4$, there is some minimal $n$ so that $h_n(c) \leq 0$ (see the 
comment 
following the proof of Lemma \ref{catchpoints}.)  If $h_n(c)=0$,  
replace $c$ with a slightly larger value, that
is, recalling our notation from the previous proof, replace $c$ with some 
number in the interval $(x_n, x_{n+1}).$  Then $h_{n+1}(c) <0$.  We stop 
once we are in the case that $h_k(c)<0$.

Also recursivly define $e_n$ by
$$e_1(c,\epsilon) = \epsilon +  \frac{4\epsilon}{c^2}
\ln(1+\frac{c}{\epsilon}),$$
$$e_n(c,\epsilon) = \epsilon +  \frac{4e_{n-1}(c,
\epsilon)}{(h_{n-1}(c))^2}
\ln(1+\frac{h_{n-1}(c)}{e_{n-1}(c, \epsilon)}).$$
The recursive definition for $e_n$ is unpleasant, but all that we shall
need is that for $c$ and $n$ fixed, $e_n(c, \epsilon) \rightarrow 0$ as
$\epsilon \rightarrow 0.$  This is easy to verify by induction.  

To begin,
we will
prove by induction that
\begin{equation}\label{claim}
x(1)-x(t) \geq \epsilon-e_n(c,\epsilon)+ 
(c-h_n(c))\sqrt{1-t}.
\end{equation}
First we show equation \eqref{claim} when $n=1$.  We have
$$x(t)-\lambda(t) \leq x(1)-\lambda (1) + c\sqrt{1-t} \leq \epsilon + c 
\sqrt{1-t}$$
which implies that
$$\dot{x}(t) \geq \frac{2}{\epsilon + c\sqrt{1-t}}.$$
Since 
$$\int_t^1 \frac{2}{a+b\sqrt{1-s}} ds = 
\frac{4}{b}\sqrt{1-t}-\frac{4a}{b^2}\ln (1+\frac{b}{a}\sqrt{1-t}),$$
integrating gives
$$x(1)-x(t) \geq \frac{4}{c}\sqrt{1-t} - \frac{4\epsilon}{c^2} 
\ln(1+\frac{c}{\epsilon}\sqrt{1-t}),$$
and so, as desired \eqref{claim} holds for $n=1$.

Next assume equation \eqref{claim} holds for $n=k$.  Then 
$$x(t) \leq x(1) -\epsilon +e_k(c,\epsilon) +(h_k(c)-c)\sqrt{1-t},$$
and so
$$x(t)-\lambda(t) \leq e_k(c,\epsilon) +h_k(c)\sqrt{1-t}.$$
This again gives us an esimate for $\dot{x}(t)$ and integrating yields
$$x(1)-x(t) \geq \frac{4}{h_k(c)}\sqrt{1-t}
-\frac{4e_k(c,\epsilon)}{h_k(c)^2}\ln
(1+\frac{h_k(c)}{e_k(c, \epsilon)}\sqrt{1-t}).$$
Thus equation \eqref{claim} holds for $n=k+1$, completing our verification of 
\eqref{claim} by induction.

Recall that $x(1) \leq c+\epsilon$.  Thus letting $t=0$ in equation 
\eqref{claim} gives
$$h_n(c)+e_n(c,\epsilon) > 0.$$
As mentioned before, by adjusting $c$ slightly if necessary, there is some 
$n$ such that $h_n(c)<0$.
Then since $e_n(c, \epsilon) \rightarrow 0$ as $\epsilon \rightarrow 
0$, there exists some $\epsilon > 0$ so that $e_n(c, \epsilon) < 
-h_n(c)$.  But this contradicts the fact that $h_n(c)+e_n(c,\epsilon) 
> 0$.
Therefore, there exists $\epsilon > 0$ so that $x(1)-\lambda(1) > \epsilon$,
for $x(t)$ the solution to \eqref{bwr} with $x_0>0$.

\end{proof}

Now we wish to run time backwards, and so we must consider the second form 
of the Loewner 
equation in the upper halfplane.  Recall that from \eqref{fw}, 
the driving term $\xi(t)$ 
generates conformal functions $f_t$, which map from $\mathbb{H}$ into 
$\mathbb{H}$.  If the image of $f_t$ is a quasislit-halfplane, 
then we 
can extend $f_t$ continuously to $\mathbb{R}$, and for each $x_0 \in 
\mathbb{R} \setminus \{\xi(0)\}$, $x(t) := f_t(x)$ is a solution to

\begin{gather}\label{fwr}
\frac{\partial}{\partial t} x(t) = \frac{-2}{x(t)- \xi(t)},\\
x(0) = x_0. \notag
\end{gather} 
Note that the solution $x(t)$ might not exist for all time.  Indeed, in 
the case that $\norm{\xi}_{\frac{1}{2}} < 4$, the following corollary shows 
that $x(t)$ will 
hit the singularity $\xi(t)$ in finite time.  We define the
hitting time 
$T(x_0)$ to be the first time that $x(t)$ equals $\xi(t)$, that is, 
$x(T(x_0))=\xi(T(x_0))$ and $x(t) \neq 
\xi(t)$ for $t < T(x_0)$.  If $x(t)$ never equals $\xi (t)$, then 
$T(x_0):=\infty$.

\begin{cor}\label{finitehit}

Let  $\xi  \in$ Lip$(\frac{1}{2})$ with  $ \norm{\xi}_{\frac{1}{2}} 
<  4$ and $\xi(0)=0$.  Suppose that $x(t)$ is a 
solution to \eqref{fwr}, with $x_0 \neq 0$.  Then $K_1x_0^2 \leq T(x_0) 
\leq K_2x_0^2$, where 
$0< K_i=K_i(\norm{\xi}_{\frac{1}{2}}) <\infty$.

\end{cor}

\begin{proof}

For $c= \norm{\xi}_{\frac{1}{2}}$, let $\epsilon = \epsilon_c >0$ be 
given as in Lemma \ref{catchpoints2}, and let $x(t)$ be the 
solution to \eqref{fwr} with $x(0)=\epsilon$.  If $T(\epsilon) >1$, 
then $\lambda(t)=\xi(1-t)-\xi(1)$ 
and $y(t) = x(1-t) - \xi(1)$ satisfy the differential 
equation 
\eqref{bwr}, with $y(0) = x(1)-\xi(1) > 0$.  
Thus Lemma 2 implies that $\epsilon = y(1)-\lambda(1) > \epsilon$.  
This is a contradiction, and so $T(\epsilon) \leq 1$.

Now suppose $x_0 >0$, with $x(t)$ again the corresponding solution to 
\eqref{fwr}.  Then by the scaling property, $\hat{\xi}(t)$ and 
$\hat{x}(t)$ satisfy equation
\eqref{fwr}, where  
$$\hat{\xi}(t):=\frac{\epsilon}{x_0}\xi(\frac{x_0^2}{\epsilon ^2}t),$$
and
$$\hat{x}(t):=\frac{\epsilon}{x_0}x(\frac{x_0^2}{\epsilon ^2}t).$$
Note that $\hat{x}(0)=\epsilon$.
Therefore $T(x_0)=\frac{x_0^2}{\epsilon ^2}T(\epsilon) \leq 
K_2x_0^2$  where $K_2=K_2(c)< \infty$.  

For the lower bound, assume first that $x_0 =1$, and assume that 
$T(1)=\delta$ 
is small.  Then since $\xi(t) \leq c\sqrt{t}$, we have $x(\delta) \leq 
c\sqrt{\delta}$.  Taking $\delta$ small enough so that 
$c\sqrt{\delta} < 
\frac{1}{2}$, let $t_0$ be the time when $x(t)=\frac{1}{2}$.  Then,
$$-\frac{1}{2}= \int_0^{t_0} \frac{-2}{x(s)-\xi(s)}ds \geq 
\frac{-2t_0}{\frac{1}{2}-c\sqrt\delta},$$
and so,
$$\frac{1}{2}(\frac{1}{2}-c\sqrt\delta)\leq 2\delta.$$ 
This leads to a contradiction if $\delta$ is sufficiently small.
Therefore $T(1) \geq K_1$ for some $K_1 = 
K_1(c) > 0.$  Then by the scaling property, $T(x_0) \geq K_1 x_0^2$.

\end{proof}

In the previous corollary, we saw that if $ \norm{\xi}_{\frac{1}{2}}
<  4$  then solutions $x(t)$ to \eqref{fwr}
will hit the singularity 
in 
finite time.  Lemma \ref{2hits} shows that there is more that is true.  For 
each finite time, there 
are 
exactly two initial points, one on each side of the singularity, so that 
the solutions to $\eqref{fwr}$ will 
hit the 
singularity at that time.

\begin{lem}\label{2hits}

Let  $\xi  \in$ Lip$(\frac{1}{2})$ with  $ \norm{\xi}_{\frac{1}{2}}
<  4$.  For each $T>0$, there exist 
exactly two real numbers $x_0, \hat{x}_0$ so that 
$x(T)=\hat{x}(T)=\xi(T)$.

\end{lem}

\begin{proof}

First notice that no two points on the same side of the singularity can 
give rise to solutions to \eqref{fwr} that will 
hit at the same time.  This follows from the fact that $\delta(t) := 
y(t)-x(t)$ is increasing in $t$ for $\xi(0)<x_0<y_0$, since
$$\dot{\delta}(t)=2\frac{y(t)-x(t)}{(y(t)-\xi(t))(x(t)-\xi(t))}.$$
Thus there are at most two points that can hit at time $T$.

Next we'll show that there is one point $x_0$ to the right of the 
singularity 
with $x(T)=\xi(T)$.  For each $n \in \mathbb{N}$, set $w_n=\xi(T) + 
\frac{1}{n}$.  Now, starting at $w_n$, run time from $T$ back to 0.  This 
corresponds to solving \eqref{bwr} with intial value $w_n$.    
Since $ \norm{\xi}_{\frac{1}{2}} <  4$, the driving term cannot catch up 
with this solution, $g_t(w_n)$, by Lemma \ref{catchpoints}, and so it is 
well-defined 
up through time $T$.  Thus, 
$x_n:=g_T(w_n)=f_T^{-1}(w_n)$ is well-defined.  Further, by Lemma 
\ref{catchpoints2}, there exists $\epsilon >0$ so that $x_n-\xi(0) > 
\epsilon$.  Therefore, $\{x_n\}$ is 
a decreasing sequence bounded below by $\xi(0)+ \epsilon$, and so it has a 
limit $x_0$.  Then $x_0 > \xi(0)$ and clearly we have $x(T)=\xi(T)$.  This 
completes the proof.

\end{proof}

\subsection*{Conformal welding with the Loewner equation}

The previous lemma allows us to define 
the welding homeomorphism $\phi: \mathbb{R} \rightarrow 
\mathbb{R}$ as the orientation-reversing map that satisfies
$\phi(x)=y$ if and only if $T(x)=T(y)$.  Thus the welding 
homeomorphism interchanges the two points which hit the singularity at the 
same time.  Note that if $\xi$ is not defined for all time, but for a finite 
interval $[0,T]$, the welding homeomorphism will not be defined on all 
$\mathbb{R}$.  However, we can overcome this technicality by setting 
$\xi(t):=\xi(T)$ for $t>T$.

This next lemma is an analogue of Lemma 3.2 found in 
\cite{RM}.

\begin{lem}\label{3.2}

Let  $\xi  \in$ Lip$(\frac{1}{2})$ with  $ \norm{\xi}_{\frac{1}{2}}
<  4$ and $\xi(0)=0$.  There exists some constant $A_0>0$, depending only 
on 
$ \norm{\xi}_{\frac{1}{2}}$, so that if 
$0 \leq x <y <z$ with $y-x=z-y$, then
\begin{equation}\label{qsymm}
\frac{1}{A_0} \leq \frac{\phi(x)-\phi(y)}{\phi(y)-\phi(z)} \leq A_0.
\end{equation}
\end{lem}

To prove this lemma, we will need the following.

\begin{lem}\label{lem}

Let $c<4$ and $0<\epsilon <1$.  Then there exists $\delta >0$ so 
that    
$$\frac{\phi(\beta)}{\phi(\alpha)} \geq 1+ \delta$$
for non-zero $\alpha$ and $\beta$ satisfying $\frac{\beta}{\alpha}
\geq 1+\epsilon$ and for any Lip$(\frac{1}{2})$ driving term $\xi$ with 
$\norm{\xi}_{\frac{1}{2}} \leq c$.

\end{lem}

\begin{proof}

Notice first that without loss of generality we can take $\alpha =-1$ and 
$\beta \leq -(1 + \epsilon)$ by the scaling property.

Suppose there is no such $\delta$ as in the statement of the lemma.  Then 
for each $n \in \mathbb{N}$ there 
exists a driving term $\xi_n$ and $\beta_n \leq -(1+\epsilon)$ so 
that 
$b_n < (1+\frac{1}{n})a_n$, where $0<a_n := \phi(-1) < b_n:= 
\phi(\beta_n)$.  
Set $T_n=T(a_n)$ and $S_n=T(b_n)$.

By Ascoli-Arzela, there exists a subsequence of $\{\xi_n\}$ which 
converges locally uniformly to $\xi$.  Note that $\xi \in$ 
Lip$(\frac{1}{2})$ with $\norm{\xi}_{\frac{1}{2}} \leq c$.  Since $T(x) 
\asymp x^2$ by Corollary \ref{finitehit}, $a_n, b_n, \beta_n, 
T_n,$ and $S_n$ are all bounded.  Hence by taking 
subsequences and renaming to avoid notational hazards, we have $a_n 
\rightarrow a, b_n \rightarrow b, \beta_n \rightarrow \beta, T_n 
\rightarrow T,$ and $S_n \rightarrow S.$  Note that $a=b$ since 
$a_n<b_n<(1+\frac{1}{n})a_n$.  If we had that $T(a)=T=T(-1)$ and 
$T(b)=S=T(\beta)$, this would give us the desired contradiction, since 
$T(-1)<T(\beta)$.  The same argument can be used to prove each of these 
four equalities, and so we will simply show that 
$T(a)=T$.  Since $\xi_n \rightarrow \xi$ locally 
uniformly, $\xi_n(T_n) \rightarrow \xi(T)$.  Hence $\lim_{n \rightarrow 
\infty} a_n(T_n) = \lim_{n \rightarrow \infty} \xi_n(T_n)= \xi(T)$, where 
$a_n(t)$ is the solution to \eqref{fwr} with $a_n(0)=a_n$.  Thus it 
remains to show that $a_n(T_n) \rightarrow a(T)$.  

Claim: Let $\epsilon >0$.  Then $a_n(T-\epsilon) \rightarrow 
a(T-\epsilon)$.

Proof of Claim:  We will assume without loss of generality that $T_n \geq 
T-\frac{\epsilon}{2}$.  
Then, $a_n(T-\epsilon)$ is well-defined and is bounded away from 
$\xi_n(T-\epsilon)$ by a factor of $\sqrt\epsilon$ by Corollary 
\ref{finitehit}.

Fix $n$ for a moment.  Then looking to solve the initial value 
problem \eqref{fwr} with the method of successive approximations, let 
$\psi_0^n \equiv a_n$ and recursively define 
$$\psi_{k+1}^n(t) = a_n + \int_0^t\frac{-2}{\psi_k^n(s)-\xi_n(s)}ds.$$
Similarly, let $\psi_k$ be the approximation for $\xi$ with initial value 
$a$.
Then for $t \in [0,T-\epsilon]$, $\psi_k^n(t) \geq a_n(t)$ 
and $\psi_k(t) \geq a(t)$.  By an easy induction 
argument, we have that for $t \in 
[0,T-\epsilon],$
$$\abs{\psi_k^n(t)-\psi_k(t)} \leq \abs{a_n-a}+(\abs{a_n-a} + 
\norm{\xi_n-\xi}_{\infty}) \sum_{j=1}^k \frac{(Bt)^j}{j!}$$
where $B$ depends only on $\epsilon$.  So, for $t \in [0,T-\epsilon]$,
\begin{equation*}
\begin{split}
\abs{a_n(t)-a(t)} &= \lim_{k \rightarrow \infty} 
 \abs{\psi_k^n(t)-\psi_k(t)} \\
 &\leq \abs{a_n-a}+(\abs{a_n-a} 
 +\norm{\xi_n-\xi}_{\infty})(\mathrm{e}^{Bt}-1). \\
\end{split}
\end{equation*}
Therefore, $a_n(T-\epsilon) \rightarrow a(T-\epsilon)$, proving the claim.

Assuming $T_n \in [T-\frac{\epsilon}{2}, T+\frac{\epsilon}{2}]$ and 
using Corollary \ref{finitehit}, we have
\begin{equation*}
\begin{split}
0 \leq a_n(T-\epsilon)-a_n(T_n) \\
 &= (a_n(T-\epsilon)-\xi_n(T-\epsilon))+(\xi_n(T-\epsilon)-\xi_n(T_n)) \\
 &\leq A\sqrt \epsilon +c\sqrt{T_n-(T-\epsilon)} \\
 &\leq A\sqrt \epsilon \\
\end{split}
\end{equation*}
where $A$ is a constant depending only on $c.$
So by the claim above,
$$0\leq a(T-\epsilon) - \lim_{n \rightarrow \infty} a_n(T_n) \leq A\sqrt 
\epsilon$$
implying that $a_n(T_n) \rightarrow a(T)$.

\end{proof}

Now we are ready for the proof of Lemma \ref{3.2}.

\begin{proof}

In this proof, $A \geq 1$ will stand for any constant which depends only 
on 
$ \norm{\xi}_{\frac{1}{2}}$.  Let $z(t)$ be the
solution to \eqref{fwr} with $z(0)=z$, and $\hat{z}(t)$ the solution to 
\eqref{fwr} with $\hat{z}(0)=\phi(z)$.  Define $x(t)$,
$y(t)$, $\hat{x}(t)$ and $\hat{y}(t)$ similarly.

First we consider the case $x=0$.  Instead of only taking $z=2y$, we 
simply 
assume that $\frac{z}{y} \in [1+\epsilon, 2]$, since we will reduce the 
next case to 
this setting.  By the scaling invariance, we
can assume that $y=1$. 
Set $T=T(1)$, and recall that $K_1 \leq T \leq K_2$ from Corollary 
\ref{finitehit}.  
Then 
$z(T)-\xi(T) \leq 2+c\sqrt K_2$.  
Abusing notation a little, we have $T(z)=T+T(z(T)-\xi(T))$, where by 
$T(z(T)-\xi(T))$ we 
mean the hitting time for the solution to \eqref{fwr} with 
initial value 
$z(T)$ and driving term $\xi(T+t)$.  By Corollary \ref{finitehit}, 
$$\phi(z)^2 \leq \frac{1}{K_1} T(\phi(z)) = \frac{1}{K_1} T(z) \leq 
\frac{K_2}{K_1}(1+(2+c\sqrt K_2)^2)$$
and similarly,
$$\phi(1)^2 \geq \frac{1}{K_2} T(\phi(1)) = \frac{1}{K_2} T(1) \geq 
\frac{K_1}{K_2}.$$
Therefore, 
$$\frac{\phi(z)}{\phi(1)} \leq A.$$
By Lemma \ref{lem}, we have 
$$\frac{\phi(z)}{\phi(1)} \geq 1+ \delta$$
where $\delta$ depends only on $c$ and $\epsilon$.  
This gives \eqref{qsymm} in the case $x=0.$

Next we consider the case where $x>0$ and $z \geq 2x$.  We will reduce 
this to case 1 by letting time run for $T=T(x)$ at which point 
$x(T)=\xi(T)$.
Since 
$$\frac{\partial}{\partial t} \log(\frac{y(t)-x(t)}{z(t)-y(t)})= 
2\frac{z(t)-x(t)}{(x(t)-\xi(t))(y(t)-\xi(t))(z(t)-\xi(t))},$$
the quotient $q(t):=\frac{y(t)-x(t)}{z(t)-y(t)}$ is increasing in $t$.  
Therefore $q(T)>1$.  Also,
$$q(T) = \frac{y(T)-x(T)}{z(T)-y(T)} \leq 
\frac{y+c\sqrt{T}}{\frac{1}{2}(z-x)} \leq
\frac{(1+c\sqrt K_2)z}{\frac{1}{4}z} \leq A.$$
Now we are back to case 1, since 
we have  $(1+\frac{1}{A})(y(T)-\xi(T)) \leq z(T)-\xi(T) \leq 
2(y(T)-\xi(T))$.  
Hence by case 1, there exists $A$ depending only on $c$, 
so 
that 
$$\frac{1}{A} \leq 
\frac{\hat{x}(T)-\hat{y}(T)}{\hat{y}(T)-\hat{z}(T)} 
\leq A.$$  
Now we would like to run time from $T$ back to 0 to give \eqref{qsymm} for 
case 2.  Since the quotient will be decreasing in $t$ as time run 
backward, we immediately get the upper bound.  For the lower bound,

$$ \frac{\phi(x)-\phi(y)}{\phi(y)-\phi(z)}
\geq \frac{\phi(x)-\phi(y)}{\hat{y}(T)-\hat{z}(T)}
\geq \frac{1}{A} \frac{\phi(x)-\phi(y)}{\hat{x}(T)-\hat{y}(T)}
\geq \frac{1}{A} \frac{\phi(x)-\phi(y)}{-\phi(y)}
\geq \frac{1}{A},$$
where Lemma \ref{lem} gives the last inequality.
Therefore \eqref{qsymm} holds for case 2.

While these first two cases required more work than in the situation in 
\cite{RM}, the final case where $x>0$ and $z-x < x$ 
follows the arguments of Lemma 3.2 in \cite{RM} without any complications.  
The idea, 
similar to the strategy used in the previous case, is to let time run 
for $S$, where $S$ is the first time that 
$x(S)-\xi(S)=z(S)-x(S)$, and to show that the quotient $q(t)$ is bounded 
on $[0,S]$.  Thus, we end up in a setting similar to 
case 2.  
It remains then to verify that case 2 still applies and to 
run time backwards from $S$ to 0, again utilizing the boundedness of 
$q(t)$.

\end{proof}

We include the statement of Lemma 2.2 from \cite{RM} below, since we will 
use it in the proof of Theorem \ref{thm:2}.  It gives a condition in terms 
of the welding homeomorphism for when a slit-halfplane is a 
quasislit-halfplane.

\begin{lem}\label{RMlem}

$\mathbb{H} \setminus \gamma[0,T]$ is a quasislit-halfplane if and only if 
there 
is a constant $1 \leq M < \infty$ such that
$$\frac{1}{M} \leq \frac{x-\xi(0)}{\xi(0)-\phi(x)} \leq M$$
for all $x> \xi(0)$ and
$$\frac{1}{M} \leq \frac{\phi(x)-\phi(y)}{\phi(y)-\phi(z)} \leq M$$
whenever $\xi(0) \leq x <y <z$ with $y-x=z-y$.  Furthermore, the quasislit 
constant $K$ of $\mathbb{H} \setminus \gamma[0,T]$ depends on $M$ only. 

\end{lem}

\subsection*{Proof of Theorem \ref{thm:2}}

\begin{proof}

By the scaling property, it suffices to show that $f_1(\mathbb{H})$ is a 
quasislit plane.  Let $n \in
\mathbb{N}$, and set $t_k = k/n$.  Following the methods in \cite{RM}, we 
wish to construct 
$\xi_n \in$ Lip$(\frac{1}{2})$ so that $\xi_n(t_k)=\xi(t_k)$ and
$\norm{\xi_n}_{\frac{1}{2}} \leq 
c:=\norm{\xi}_{\frac{1}{2}}$.  There are at least two ways to proceed.  
The first is by linear interpolation, and this is the method we will use.  
Alternatively, setting 
$c_k=(\xi(t_k)-\xi(t_{k+1}))\sqrt n$, we can define $\hat\xi_n(t)$ for 
$t \in [0,1]$ by $\hat\xi_n|_{[t_k,t_{k+1}]}(t) = c_k \sqrt{t_{k+1}-t} + 
\xi(t_{k+1})$.  Although $\hat\xi_n \in$ 
Lip$(\frac{1}{2})$, it may not be true that 
$\norm{\hat\xi_n}_{\frac{1}{2}} \leq
c$.  However, it is possible to complete the proof using this construction 
for 
$\hat\xi_n$ by considering the larger space of locally Lip$(\frac{1}{2})$ 
functions and verifying that all the lemmas remain true for these 
functions as well.  The benefit to using this construction is that we know 
slightly more about the generated domains.   If $\hat\phi_t^k$ is the 
map generated by $\hat\xi_n(t_k+t)=
c_k\sqrt{\frac{1}{n}-t}+\alpha_{k+1}$ for $t \in [0,
\frac{1}{n}]$, then $\hat\phi_{\frac{1}{n}}^k$ is a map from
$\mathbb{H}$ onto the upper halfplane slit by a line segment whose
angle with the real line is bounded away from 0 and $\pi$. 

Using our first method of linear interpolation, we set 
$m_k=n(\xi(t_{k+1})-\xi(t_k))$ and define $\xi_n(t)$ for
$t \in [0,1]$ by $\xi_n|_{[t_k,t_{k+1}]}(t) =m_k(t-t_k)+\xi(t_k)$.  First 
we check that 
$\norm{\xi_n}_{\frac{1}{2}} \leq c$.  Let $x,y 
\in [0,1]$.  If $x,y \in [t_k, t_k+1]$ for some $k$, then clearly
$\abs{\xi_n(y)-\xi_n(x)} \leq c \sqrt{\abs{y-x}}$.  So assume that $t_j 
\leq x \leq t_{j+1} \leq t_k  \leq y \leq t_{k+1}$, and assume without 
loss of 
generality that $\xi_n(y) \geq \xi_n(x)$.  If we maximize the function 
$h(x,y):=\xi_n(y) - \xi_n(x) - c\sqrt{y-x}$ over $(x,y) \in [t_j, 
t_{j+1}] \times [t_k,t_{k+1}],$ we find that $h(x,y) \leq 0$, as desired.

Let $\phi_t^k$ be the maps generated by $\xi_n(t_k+t)= 
m_kt+\xi(t_{k})$ for $t \in [0, 
\frac{1}{n}]$.  Then $\phi^k:=\phi_{\frac{1}{n}}^k$ is a map from 
$\mathbb{H}$ onto the upper halfplane slit by a smooth curve which makes 
an 
angle of $\frac{\pi}{2}$ with the real line.  If 
$f_t^n$ is the map generated by $\xi_n$ for $t \in [0,1]$, we have that 
$f_1^n=\phi^n \circ \phi^{n-1} \circ \cdots \circ \phi^2 \circ \phi^1$.
Hence, $f_1^n(\mathbb{H})$ is a slit-halfplane.  By Corollary 
\ref{finitehit}, the first condition of Lemma \ref{RMlem} is satisfied, 
while the second condition is a result of
Lemma \ref{3.2}.  Therefore, we have that $f_1^n(\mathbb{H})$ is a 
$K$-quasislit-halfplane, with $K$ independent of $n$.  By compactness of 
the space of $K$-quasislit-halfplanes, we have that $f_1(\mathbb{H})$ is a 
quasislit-halfplane.

\end{proof}

\bibliographystyle{amsplain}
\bibliography{bib}
\end{document}